\documentclass[a4paper,12pt]{article}
\usepackage{latexsym}
\usepackage{amsfonts}
\newtheorem{theorem}{Theorem}[section]
\newtheorem{corollary}[theorem]{Corollary}
\newtheorem{definition}{Definition}

\newtheorem{observation}[theorem]{Observation}
\newtheorem{proposition}[theorem]{Proposition}

\title{An answer to a question of Coleman on scattered
sets\thanks{1991 Math.\ Subject Classification ---
Primary: 54G12, 54H05; Secondary: 68U05, 68U10. \protect\newline Key
words and phrases --- scattered, nowhere dense, semi-$T_D$-space,
density topology, digital topology, fenestration, trace space.
\protect\newline Research partially supported by the Ella and Georg
Ehrnrooth Foundation at Merita Bank, Finland.}}
\author{Julian Dontchev\\Department of Mathematics\\University
of Helsinki \\Yliopistonkatu 5\\00014 Helsinki\\Finland \and
Maximilian Ganster\\Department of Mathematics\\Graz University
of Technology\\Steyrergasse 30\\A-8010 Graz\\Austria}
\date{}
\begin{document}
\baselineskip=20pt plus 1pt minus 1pt
\newcommand{\fxy}{$f \colon (X,\tau) \rightarrow (Y,\sigma)$}
\maketitle
\begin{abstract}
The aim of this paper is to show that every scattered subset of
a dense-in-itself semi-$T_D$-space is nowhere dense. We are thus
able to answer a recent question of Coleman \cite{C1} in the
affirmative. In terms of Digital Topology, we prove that in
semi-$T_D$-spaces with no open screen, trace spaces have no
consolidations.
\end{abstract}

\section{Introduction}\label{s1}

It is well-known that in dense-in-themselves $T_1$-spaces, all
scattered subsets are nowhere dense. This result was established
by Kuratowski in the proof that in $T_1$-spaces the finite union
of scattered subsets is scattered.

In a recent paper Coleman asked the following question
\cite[Question 4]{C1}: Is it true that in dense-in-themselves,
$T_D$-spaces all scattered sets are nowhere dense? In what
follows, we will show that even in dense-in-themselves
semi-$T_D$-spaces all $\alpha$-scattered sets are nowhere dense.

The question of Coleman is in fact very well motivated not only
because it is interesting to know how low one can go on the
separations below $T_1$ and still have the scattered sets being
nowhere dense but also from a `digital point of view'. In terms
of Digital Topology, we will prove that in semi-$T_D$-spaces with
empty open screen, trace spaces have no consolidations.

In Digital Topology several spaces that fail to be $T_1$ are very
often important in the study of the geometric and topological
properties of digital images \cite{KR1,K1}. Such is in fact the
case with the major building block of the digital n-space -- the
{\em digital line} or the so called {\em Khalimsky line}. This
is the set of the integers, $\mathbb Z$, equipped with the
topology $\cal K$, generated by ${\cal G}_{\cal K} = \{ \{ 2n-1,
2n, 2n+1 \} \colon n \in {\mathbb Z} \}$.

A {\em fenestration} \cite{K1} of a space $X$ is a collection of
disjoint nonempty open sets whose union is dense. The {\em
consolidation} $A^+$ \cite{K1} of a set $A$ is the interior of
its closure. When there is a fenestration of a space $(X,\tau)$
by singletons, the space $(X,\tau)$ is called {\em
$\alpha$-scattered} \cite{DGR1} or a {\em trace space} \cite{K1}.
For example, in the digital line the collection $\{ \{ n \} :  
n \in {\mathbb Z}$ and $n$ is odd$\}$ is a fenestration of
$({\mathbb Z},{\cal K})$. All scattered sets are
$\alpha$-scattered by not vice versa \cite{DGR1}. In $T_0$-spaces
without isolated points, we may encounter a trace space which
fails to be nowhere dense \cite{C1}. Nevertheless, as we will
show, with the presence of the very weak separation 'semi-$T_D$',
in spaces with no isolated points we have all $\alpha$-scattered
sets being nowhere dense.

A topological space $X$ is a called a {\em $T_{D}$-space} if
every singleton is locally closed or equivalently if the derived
set $d(x)$ is closed for every $x \in X$. Recall that $X$ is a
{\em semi-$T_{D}$-space} if every singleton is open or
semi-closed \cite{D1}. Recall that a subset $A$ of a space
$(X,\tau)$ is called {\em locally dense} \cite{CM1} if $A
\subseteq A^+$. Note that every open and every dense set is
locally dense.

\section{When is $\cal N$ finer than $\cal S$?}\label{s2}

Recall that a topological ideal $\cal I$ is a nonempty collection
of sets of a space $(X,\tau)$ closed under heredity and finite
additivity. For example, the families $\cal N$ (of all nowhere
dense sets) and $\cal F$ (of all finite sets) always form ideals
while the family $\cal S$ of all scattered sets is an ideal if
and only if the space is $T_0$.

\begin{proposition}\label{p1}
For a topological space $(X,\tau)$ the following conditions are
equivalent:

{\rm (1)} $X$ is a dense-in-itself semi-$T_D$-space.

{\rm (2)} Every singleton is nowhere dense.

{\rm (3)} ${\cal F} \subseteq {\cal N}$.

{\rm (4)} There are no locally dense singletons in $X$.
\end{proposition}

{\em Proof.} (1) $\Rightarrow$ (2) Let $x \in X$. Since $X$ is
a semi-$T_D$-space, $\{ x \}$ is open or semi-closed \cite{D1}.
Since $X$ is dense-in-self, $\{ x \}$ is semi-closed. On the
other hand, in any topological space every singleton is locally
dense (= preopen) or nowhere dense \cite{JR1}. If $\{ x \}$ is
preopen, then it must be (due to semi-closedness) regular open.
As $X$ has no isolated points, we conclude that $\{ x \}$ is
nowhere dense.

(2) $\Rightarrow$ (3) Obvious, since the ideal of nowhere dense
sets is closed under finite additivity.

(3) $\Rightarrow$ (4) Follows easily from the fact that
singletons are either locally dense or nowhere dense.

(4) $\Rightarrow$ (1) If some point $x \in X$ were isolated, then
it would be locally dense. This shows that $X$ is
dense-in-itself. That $X$ is a semi-$T_D$-space follows easily
from the fact that nowhere dense sets are semi-closed.

Recall that a subset $A$ of a topological space $(X,\tau)$ is
called {\em $\beta$-open} if $A$ is dense in some regular closed
subspace of $X$. Note that every locally dense set is
$\beta$-open.

\begin{observation}\label{p2}
{\rm (i)} Every $\beta$-open subset of a dense-in-itself
semi-$T_D$-space is also dense-in-itself and semi-$T_D$.

{\rm (ii)} Let $(X_i,\tau_i)_{i \in I}$ be a family of
topological spaces such that at least one of them is a
dense-in-itself semi-$T_D$-space. Then the product space $X =
\prod_{i \in I} X_i$ is also dense-in-itself and semi-$T_D$.
\end{observation}

\begin{theorem}\label{t1}
If a topological space $(X,\tau)$ is dense-in-itself and
semi-$T_D$, then every $\alpha$-scattered subset of $X$ is
nowhere dense.
\end{theorem}

{\em Proof.} Let $A \subseteq X$ be $\alpha$-scattered. Assume
that $A^+$ is nonempty, i.e., there exists a nonempty $U \in
\tau$ such that $U \subseteq {\rm cl} (A)$. Since $(A,\tau|A)$
is $\alpha$-scattered, $U$ meets $I(A)$, the set of all isolated
points of $(A,\tau|A)$. Let $x \in U \cap I(A)$ and let $V$ be
an open subset of $(X,\tau)$ such that $V \cap A = \{ x \}$. Set
$W = U \cap V$. Note that $W \subseteq V \cap \overline{A}
\subseteq \overline{V \cap A}  = \overline{\{ x \}}$ and so $\{
x \}$ has nonempty consolidation, i.e. it is not nowhere dense
in $X$. By Proposition~\ref{p1}, we have a contradiction. Hence,
$A$ is nowhere dense. $\Box$

The digital interpretation of Theorem~\ref{t1} is as follows: In
semi-$T_D$-spaces with no isolated points, i.e.,
semi-$T_D$-spaces with empty open screens, the trace spaces have
empty consolidations.

Now, we can apply the result above in order to show that the
$\alpha$-scattered subsets of the density topology are in fact
its Lebesgue null set.

\begin{definition}\label{d2}
{\em A measurable set $E \subseteq {\mathbb R}$ has density $d$
at $x \in {\mathbb R}$ if $$\lim_{h \rightarrow 0} \frac{m(E \cap
[x-h,x+h])}{2h}$$ exists and is equal to $d$. Set $\phi(E) = \{
x \in {\mathbb R} \colon d(x,E) = 1 \}$. The open sets of the
density topology $\cal T$ are those measurable sets $E$ that
satisfy $E \subseteq \phi(E)$. Note that the density topology
$\cal T$ is finer than the usual topology on the real line.}
\end{definition}

\begin{corollary}\label{c1}
The trace spaces (i.e., the $\alpha$-scattered subsets) of the
density topology are precisely its Lebesgue null set.
\end{corollary}

{\em Proof.} Follows from Theorem~\ref{t1} and the fact that a
subset $A$ of the density topology is nowhere dense if and only
if it is a Lebesgue null set \cite{T1}. $\Box$

\
e-mail: {\tt dontchev@cc.helsinki.fi}, {\tt
ganster@weyl.math.tu-graz.ac.at}
\
\end{document}